\documentclass[12pt]{amsart}
\usepackage{amssymb,amsmath,amscd}
\makeatletter
\def\@cite#1#2{{\m@th\upshape\bfseries%
[{#1\if@tempswa{\m@th\upshape\mdseries, #2}\fi}]}}
\makeatother
\numberwithin{equation}{section}

\newtheorem{thm}{Theorem}[section]
\newtheorem{lem}[thm]{Lemma}

\newtheorem{prop}[thm]{Proposition}
\newcommand{\Prf}{\noindent\textbf{Proof.\ }}
\newcommand{\pg}{\mathcal{C}_\Gamma}

\begin{document}
\title[C*-algebras generated by composition operators]{$C^*$-algebras generated by groups of composition operators}
\author{Michael~T.~Jury} 
\address{Department of Mathematics\\
         University of Florida\\
         Gainesville, Florida 32603}
\email{mjury@math.ufl.edu}
\date{26 Sep 2005}
\thanks{*Research partially supported by NSF VIGRE grant DMS-9983601}
\begin{abstract}
We compute the C*-algebra generated by a group of composition operators
acting on certain reproducing kernel Hilbert spaces over the disk, where the
symbols belong to a non-elementary Fuchsian group.  We show that such a
C*-algebra contains the compact operperators, and its quotient is
isomorphic to the crossed product C*-algebra determined by the action
of the group on the boundary circle.  In addition we show that the
C*-algebras obtained from composition operators acting on a natural
family of Hilbert spaces are in fact isomorphic, and also determine
the same \emph{Ext}-class, which can be related to known extensions of
the crossed product.
\end{abstract}
\thanks{2000 {\it Mathematics Subject Classification.} 20H10, 46L55, 46L80, 47B33, 47L80}
\maketitle

\section{Introduction}

The purpose of this paper is to begin a line of investigation
suggested by recent work of Moorhouse et al.  \cite{moorhouse-thesis,
  jm-tk-preprint, jm-tk-bm-preprint}: to describe, in as much detail
as possible, the $C^*$-algebra generated by a set of composition
operators acting on a Hilbert function space.  In this paper we
consider a class of examples which, while likely the simplest cases
from the point of view of composition operators, nonetheless produces
$C^*$-algebras which are of great interest both intrinsically and for
applications.  In fact, the $C^*$-algebras we obtain are objects of
current interest among operator algebraists, with applications in the
study of hyperbolic dynamics \cite{emerson-thesis}, noncommutative
geometry \cite{MR88k:58149}, and even number theory \cite{MR1931172}.

Let $f\in H^2(\mathbb{D})$.  For an analytic function $\gamma:\mathbb{D}\to\mathbb{D}$,
the \emph{composition operator with symbol $\gamma$} is the linear
operator defined by
$$
(C_\gamma f) (z)=f(\gamma(z))
$$ 
In this paper, we will be concerned with the $C^*$-algebra
$$
\pg=C^*(\{C_\gamma : \gamma\in\Gamma\})
$$
where $\Gamma$ is a discrete group of (analytic) automorphisms of
$\mathbb{D}$ (i.e. a \emph{Fuchsian group}).  For reasons to be described
shortly, we will further restrict ourselves to \emph{non-elementary}
Fuchsian groups (i.e. groups $\Gamma$ for which the $\Gamma$-orbit of
$0$ in $\mathbb{D}$ accumulates at at least three points of the unit circle
$\mathbb{T}$.)  Our main theorem shows that $\pg$ contains the compact
operators, and computes the quotient $\pg/\mathcal{K}$:

\begin{thm}  
  Let $\Gamma$ be a non-elementary Fuchsian group, and let $\pg$
  denote the $C^*$-algebra generated by the set of composition
  operators on $H^2$ with symbols in $\Gamma$.  Then there is an exact
  sequence
\begin{equation}\begin{CD}
0 @>>> \mathcal{K} @>\iota>> \pg @>\pi >> C(\mathbb{T})\times\Gamma @>>> 0
\end{CD}
\end{equation}
\end{thm}
\noindent 
Here $C(\mathbb{T})\times\Gamma$ is the crossed product
$C^*$-algebra obtained from the action $\alpha$ of $\Gamma$ on
$C(\mathbb{T})$ given by
$$
\alpha_\gamma(f)(z)=f(\gamma^{-1}(z))
$$ (Since the action of $\Gamma$ on $\partial\mathbb{D}$ is amenable,
the full and reduced crossed products coincide; we will discuss this
further shortly.)  We will recall the relevant definitions and facts we require
in the next section.  

There is a similar result for the $C^*$-algebras 
$$
\pg^n =C^*\left(\{C_\gamma \in \mathcal{B}(A^2_n)|\gamma\in\Gamma\}\right),
$$  
acting on the family of reproducing kernel Hilbert spaces $A^2_n$ (defined below), 
namely there is an extension
\begin{equation}\begin{CD}
0 @>>> \mathcal{K} @>>> \pg^n @>>> C(\mathbb{T})\times\Gamma @>>> 0
\end{CD}
\end{equation}
and we will show that each of these extensions represents the same element of the
$Ext$ group
$Ext(C(\partial\mathbb{D})\times\Gamma,\mathcal{K})$, and we will also
prove that $\pg$ and $\pg^n$ are isomorphic as $C^*$-algebras.  In fact
we obtain a stronger isomorphism result, namely that \emph{any}
unital $C^*$-extension of $C(\partial\mathbb{D})\times\Gamma$ which
defines the same $Ext$-class as $\pg$ is isomorphic to $\pg$:
\begin{thm}
Let $x\in \rm{Ext}(C(\partial\mathbb{D})\times\Gamma)$ denote the class of the extension
\begin{equation}\nonumber\begin{CD}
0 @>>> \mathcal{K} @>>> \pg @>>> C(\partial\mathbb{D})\times\Gamma @>>> 0.
\end{CD}
\end{equation}
If $e\in \rm{Ext}(C(\partial\mathbb{D})\times\Gamma)$ is a unital extension represented by 
\begin{equation}\nonumber\begin{CD}
0 @>>> \mathcal{K} @>>> E  @>>> C(\partial\mathbb{D})\times\Gamma @>>> 0
\end{CD}
\end{equation}
and $e=x$, then $E\cong \pg$ as $C^*$-algebras.
\end{thm}

Finally, we will compare the extension determined by $\pg$ to tow
other recent constructions of extensions of
$C(\partial\mathbb{D})\times\Gamma$.  We show that the $Ext$-class of
$\pg$ coincides with the class of the $\Gamma$-equivariant Toeplitz
extension of $C(\partial\mathbb{D})$ constructed by J. Lott
\cite{lott-preprint}, and differs from the extension of crossed
products by cocompact groups constructed by H. Emerson \cite{emerson-thesis}.  Finally
we show that this extension in fact gives rise to a
$\Gamma$-equivariant $KK_1$-cycle for $C(\partial\mathbb{D})$ which
also accords with the construction in \cite{lott-preprint}.

\section{Preliminaries}

We will consider C*-algebras generated by composition operators which act on a family of reproducing kernel Hilbert spaces on the unit disk.  Specifically we will consider the spaces of analytic functions $A^2_n$, where $A^2_n$ is the space with reproducing kernel
$$
k_n(z,w)=(1-z\overline{w})^n
$$
When $n=1$ this space is the Hardy space $H^2$, and its norm is given by
$$
\|f\|^2 =\lim_{r\to 1}\frac{1}{2\pi}\int_0^{2\pi} |f(re^{i\theta})|^2\, d\theta
$$
For $n\geq 2$, the norm on $A^2_n$ is equivalent to 
$$
\|f\|^2=\frac{1}{\pi}\int_\mathbb{D}|f(z)|^2 (1-|z|^2)^{n-2}\, dA(z)
$$
though the norms are equal only for $n=2$ (the Bergman space).  

An analytic function $\gamma:\mathbb{D}\to\mathbb{D}$ defines a
\emph{composition operator} $C_\gamma$ on $A^2_n$ by
$$
(C_\gamma f)(z)=f(\gamma(z))
$$
In this paper, we will only consider cases where $\gamma$ is a
M\"{o}bius transformation; in these cases $C_\gamma$ is easily seen to
be bounded, by changing variables in the integrals defining the
norms.  An elementary calculation shows that if $\gamma :\mathbb{D}\to \mathbb{D}$ is analytic, then
$$
C_\gamma^* k_w(z) = k_{\gamma(w)}(z)
$$
where $k$ is any of the reproducing kernels $k_n$.

We recall here the definitions of the full and reduced crossed product
$C^*$-algebras; we refer to \cite{MR81e:46037} for details.  Let a
group $\Gamma$ act by homeomorphisms on a compact Hausdorff space $X$.
This induces an action of $\Gamma$ on the commutative $C^*$-algebra
$C(X)$ via
$$
(\gamma\cdot f)(x)=f(\gamma^{-1}\cdot x)
$$  
The \emph{algebraic crossed product} $C(X)\times_{alg} \Gamma$ consists
of formal finite sums $f=\sum_{\gamma\in\Gamma}f_\gamma [\gamma]$, where
$f_\gamma\in C(\mathbb{T})$ and the $[\gamma]$ are formal symbols.  Multiplication is defined in $C(X)\times_{alg}
\Gamma$ by
$$
\left(\sum_{\gamma\in\Gamma}f_\gamma [\gamma] \right)
\left(\sum_{\gamma\in\Gamma}f_\gamma^\prime [\gamma^\prime] \right) =
\sum_{\delta\in\Gamma}\sum_{\gamma\gamma^\prime=\delta} f_\gamma
(\gamma\cdot f^\prime_{\gamma^\prime})[\delta]
$$
For $f=\sum_\gamma f_\gamma [\gamma]$, define $f^*\in
C(X)\times_{alg} \Gamma$ by
$$
f^* = \sum_{\gamma\in\Gamma} (\gamma \cdot \overline{f_{\gamma^{-1}}})
[\gamma]
$$
With this multiplication and involution, $C(X)\times_{alg} \Gamma$
becomes a $*$-algebra, and we may construct a $C^*$-algebra by closing
the algebraic crossed product with respect to a $C^*$-norm.

To obtain a $C^*$-norm, one constructs $*$-representations of
$C(X)\times_{alg}\Gamma$ on Hilbert space.  To do this, we first fix a faithful representation $\pi$ of $C(X)$ on a Hilbert space $\mathcal{H}$.  We then construct a representation $\sigma$ of the algebraic crossed product on $\mathcal{H}\otimes\ell^2(\Gamma)=\ell^2(\Gamma, \mathcal{H})$ as follows:  define a representation $\tilde{\pi}$ of $C(X)$ by its action on vectors $\xi\in\ell^2(\Gamma, \mathcal{H})$
$$
(\tilde{\pi}(f))(\xi)(\gamma)=\pi(f\circ\gamma)\xi(\gamma)
$$
Represent $\Gamma$ on $\ell^2(\Gamma,\mathcal{H})$ by left translation:
$$
(U(\gamma))(\xi)(\eta)=\xi(\gamma^{-1}\eta)
$$
The representation $\sigma$ is then given by 
$$
\sigma\left(\sum f_\gamma [\gamma]\right)=\sum \tilde{\pi}(f_\gamma) U(\gamma)
$$
The closure of $C(X)\times_{alg}\Gamma$ with respect the norm induced by this representation 
is called the \emph{reduced crossed product} of $\Gamma$ and $C(X)$,
and is denoted $C(X)\times_r\Gamma$.  The \emph{full crossed product},
denote $C(X)\times \Gamma$, is obtained by taking the closure of the
algebraic crossed product with respect to the maximal $C^*$-norm
$$
\|f\|=\sup_\pi \|\pi(f)\|
$$
where the supremum is taken over \emph{all} $*$-representations
$\pi$ of $C(X)\times_{alg}\Gamma$ on Hilbert space.  When $\Gamma$ is
discrete, $C(X)\times \Gamma$ contains a canonical subalgebra
isomorphic to $C(X)$, and there is a natural surjective
$*$-homomorphism $\rho:C(X)\times\Gamma \to C(X)\times_r \Gamma$.

The full crossed product is important because of its universality with respect to covariant representations.  A \emph{covariant representation} of the pair $(\Gamma, X)$ consists of a faithful representation $\pi$ of $C(X)$ on Hilbert space, together with a unitary representation $u$ of $\Gamma$ on the same space satisfying the covariance condition
$$
 u(\gamma)\pi(f)u(\gamma)^* = f\circ \gamma^{-1}
$$
for all $f\in C(X)$ and all $\gamma\in\Gamma$.  If $\mathcal{A}$ denotes the C*-algebra generated by the images of $\pi$ and $u$, then there is a surjective *-homomorphism from $C(X)\times\Gamma$ to $\mathcal{A}$; equivalently, any C*-algebra generated by a covariant representation is isomorphic to a quotient of the full crossed product.

We now collect properties of group actions on topological spaces which
we will require in what follows.  Let a group $\Gamma$ act by
homeomorphisms on a locally compact Hausdorff space $X$; The action of
$\Gamma$ on $X$ is called \emph{minimal} if the set $\{\gamma\cdot x |
\gamma\in \Gamma\}$ is dense in $X$ for each $x\in X$, and called
\emph{topologically free} if, for each $\gamma\in \Gamma$, the set of
points fixed by $\gamma$ has empty interior.

Suppose now $\Gamma$ is discrete and $X$ is compact.  Let
$\text{Prob}(\Gamma)$ denote the set of finitely supported probability
measures on $\Gamma$.  We say $\Gamma$ \emph{acts amenably on} $X$ if
there exists a sequence of weak-* continuous maps $b_x^n:X\to
\text{Prob}(\Gamma)$ such that for every $\gamma\in \Gamma$,
$$
\lim_{n\to \infty} \sup_{x\in X} \|\gamma\cdot b_x^n - b_{\gamma\cdot x}^n\|_1
=0
$$
where $\Gamma$ acts on the functions $b_x^n$ via $(\gamma\cdot
b_x^n)(z)=b_x^n(\gamma^{-1}\cdot z)$, and $\|\cdot \|_1$ denotes the
$l^1$-norm on $\Gamma$.

\begin{thm}\label{T:as}\cite{MR94m:46101}
  Let $\Gamma$ be a discrete group acting on a compact Hausdorff space
  $X$, and suppose that the action is topologically free.  If $\mathcal{J}$ is
  an ideal in $C(X)\times \Gamma$ such that $C(X)\cap \mathcal{J}=0$, then
  $\mathcal{J}\subseteq \mathcal{J}_\lambda$, where $\mathcal{J}_\lambda$ is the kernel of the
  projection of the full crossed product onto the reduced crossed
  product.
\end{thm}

\begin{thm}\label{T:amenability}\cite{MR2000g:19004, MR86j:22014}
  If the $\Gamma$ acts amenably on $X$, then the full and reduced
  crossed product $C^*$-algebras coincide, and this crossed product is
  nuclear.  
\end{thm}

In this paper we are concerned with the case $X=\mathbb{T}$ and $G=\Gamma$,
a non-elementary Fuchsian group.  The action of $\Gamma$ on $\mathbb{T}$ is
always amenable, so by Theorem~\ref{T:amenability} the full and
reduced crossed products coincide, $C(\partial\mathbb{D})\times\Gamma$ is nuclear.  

We also require some of the basic terminology concerning Fuchsian groups.  Fix a base point $z_0\in\mathbb{D}$.  The \emph{limit set} of $\Gamma$ is the set accumulation points of the orbit $\{\gamma(z_0) : \gamma\in\Gamma\}$ on the boundary circle; it is closed (and does not depend on the choice of $z_0$.  The limit set can be one of three types:  it is either finite (in which case it consists of at most two points), a totally disconnected perfect set (hence uncountable), or all of the circle.  In the latter case we say $\Gamma$ is of the \emph{first kind}, and of the \emph{second kind} otherwise.  If the limit set is finite, $\Gamma$ is called \emph{elementary}, and \emph{non-elementary} otherwise.   

Though we do not require it in this paper, it is worth recording that
if $\Gamma$ is of the first kind, then
$C(\partial\mathbb{D})\times\Gamma$ is simple, nuclear, purely
infinite, and belongs to the bootstrap category $\mathcal{N}$ (and is clearly
separable and unital).  Hence, it satisfies the hypotheses of the
Kirchberg-Phillips classification theorem, and is classified up to
isomorphism by its (unital) $K$-theory \cite{MR98a:46085,
MR2000k:46093}.

Finally, we recall briefly the basic facts about extensions of C*-algebras and the \emph{Ext} functor.  An exact sequence of C*-algebras 
\begin{equation}\begin{CD}
0 @>>> B @>>> E @>>> A @>>> 0
\end{CD}
\end{equation}
is called an \emph{extension} of $A$ by $B$.  (In this paper we consider only extensions for which $B=\mathcal{K}$, the C*-algebra of compact operators on a separable Hilbert space.)  Associated to an extension of $A$ by $\mathcal{K}$ is a $*$-homomorphism $\tau$ from $A$ to the Calkin algebra $\mathcal{Q}(\mathcal{H})=\mathcal{B}(\mathcal{H})/\mathcal{K}(\mathcal{H})$; this $\tau$ is called the \emph{Busby map} associated to the extension.  Conversely, given a map $\tau:A\to \mathcal{Q}(\mathcal{H})$ there is a unique extension having Busby map $\tau$; we will thus speak of an extension and its Busby map interchangeably.  Two extensions of $A$ by $\mathcal{K}$ with Busby maps $\tau_1:A\to \mathcal{Q}(\mathcal{H}_1)$ and $\tau_2:A\to \mathcal{Q}(\mathcal{H}_2)$ are \emph{strongly unitarily equivalent} if there is a unitary $u:\mathcal{H}_1\to \mathcal{H}_2$ such that 
\begin{equation}\label{E:ue}
\pi(u)\tau_1(a)\pi(u)^*=\tau_2(a)
\end{equation}
for all $a\in A$ (here $\pi$ denotes the quotient map from $\mathcal{B}(\mathcal{H}_1)$ to $\mathcal{Q}(\mathcal{H}_1)$).  We say $\tau_1$ and $\tau_2$ are \emph{unitarily equivalent}, written $\tau_1\sim_u \tau_2$, if (\ref{E:ue}) holds with $u$ replaced by some $v$ such that $\pi(v)$ is unitary.  An extension $\tau$ is \emph{trivial} if it lifts to a $*$-homomorphism, that is there exists a $*$-homomorphism $\rho:A\to \mathcal{B}(\mathcal{H})$ such that $\tau(a)=\pi(\rho(a))$ for all $a\in A$.  Two extensions $\tau_1, \tau_2$ are \emph{stably equivalent} if there exist trivial extensions $\sigma_1, \sigma_2$ such that $\tau_1\oplus\sigma_1$ and $\tau_2\oplus \sigma_2$ are strongly unitarily equivalent; stable equivalence is an equivalence relation.  If $A$ is separable and nuclear (which will always be the case in this paper) then the stable equivalence classes of extensions of $A$ by $\mathcal{K}$ form an abelian group (where addition is given by direct sum of Busby maps) called $Ext(A,\mathcal{K})$, abbreviated $Ext(A)$.  Finally, each element of $Ext(A)$ determines an \emph{index homomorphism} $\partial:K_1(A)\to \mathbb{Z}$ obtained by lifting a unitary in $M_n(A)$ representing the $K_1$ class to a Fredholm operator in $M_n(E)$ and taking its Fredholm index.


\section{The Extensions $\pg$ and $\pg^n$}

\subsection{The Hardy space}
\noindent This section is devoted to the proof of our main theorem, in the case of the Hardy space:
\begin{thm}\label{T:main}  
  Let $\Gamma$ be a non-elementary Fuchsian group, and let $\pg$
  denote the $C^*$-algebra generated by the set of composition
  operators on $H^2$ with symbols in $\Gamma$.  Then there is an exact
  sequence
\begin{equation}\begin{CD}
0 @>>> \mathcal{K} @>\iota>> \pg @>\pi >> C(\mathbb{T})\times \Gamma @>>> 0
\end{CD}
\end{equation}

\end{thm}
\vskip.2in

While the proof in the general case of $A^2_n$ follows similar lines,
we prove the $H^2$ case first since it is technically simpler and
illustrates the main ideas.  The proof splits into three parts: first,
we prove that $\pg$ contains the unilateral shift $S$ and hence the
compacts (Proposition~\ref{P:shift}).  We then prove that the quotient
$\pg/\mathcal{K}$ is generated by a covariant representation of the topological
dynamical system $(C(\mathbb{T}), \Gamma)$.  Finally we prove that the
$C^*$-algebra generated by this representation is all of the crossed
product $C(\mathbb{T})\times \Gamma$.

We first require two computational lemmas.

\begin{lem}\label{L:factor}
Let $\gamma$ be an automorphism of $\mathbb{D}$ with $a=\gamma^{-1}(0)$ and
let 
$$
f(z)=\frac{1-\overline{a}z}{(1-|a|^2)^{1/2}}
$$
Then
$$
C_\gamma C_\gamma^* = M_f M_f^*
$$
where $M_f$ denotes the operator of multiplication by $f$.
\end{lem}
\Prf It suffices to show
\begin{equation}\label{E:bilin}
\langle C_\gamma C_\gamma^* k_w, k_z\rangle = \langle M_f M_f^* k_w,
k_z\rangle
\end{equation}
for all $z,w$ in $\mathbb{D}$.  Since $C_\gamma^*k_\lambda =k_{\gamma(\lambda)}$
and $M_f^* k_\lambda = \overline{f(\lambda)}k_\lambda$, (\ref{E:bilin})
reduces to the well-known identity
\begin{equation}\label{E:phi_ident}
\frac{1}{1-\gamma(z)\overline{\gamma(w)}}=\frac{(1-\overline{a}z)(1-a\overline{w})}{(1-|a|^2)(1-z\overline{w})}.
\end{equation}

\begin{lem}\label{L:commute}
Let $\gamma(z)$ be an automorphism of $\mathbb{D}$ and let $S=M_z$.  Then
$$
M_{\gamma} C_\gamma = C_\gamma S
$$
\end{lem}
\Prf For all $g\in H^2$, 
\begin{align}
\nonumber (C_\gamma S)(g)(z) &= \gamma(z)g(\gamma(z)) \\
\nonumber                  &= (M_\gamma C_\gamma)(g)(z)
\end{align}

\begin{prop}\label{P:shift}
  Let $\Gamma$ be a non-elementary Fuchsian group.  Then $\pg$
  contains the unilateral shift $S$.

\end{prop}
\Prf Let $\Lambda$ denote the limit set of $\Gamma$.  Since $\Gamma$
is non-elementary, there exist three distinct points $\lambda_1$,
$\lambda_2$, $\lambda_3\in\Lambda$.  For each $\lambda_i$ there exists a sequence $\gamma_{n,i}$ such that $\gamma_{n,i}(0)\to \lambda_i$.   Let $a_i^n
=\gamma_{n,i} (0)$.  By Lemma~\ref{L:factor},
$$
(1-|a_i^n|^2)C_{\gamma^{-1}_{n,i}}C_{\gamma^{-1}_{n,i}}^* = (1-\overline{a_i^n} S)(1-\overline{a_i^n} S)^*
$$
As $n\to \infty$, $a_i^n \to \lambda_i$ and the right-hand side
converges to 
$$
1-\overline{\lambda}S-\lambda S^* +SS^*
$$
in $\pg$.  Taking differences of these operators
for different values of $i$, we see that $\Re [\mu_1 S],\ \Re [\mu_2
S]\in \pg$, with $\mu_1 =\overline{\lambda_1 -\lambda_2}$ and
$\mu_2=\overline{\lambda_2 -\lambda_3}$.  Note that since $\lambda_1$,
$\lambda_2$, $\lambda_3$ are distinct points on the circle, the
complex numbers $\mu_1$ and $\mu_2$ are linearly independent over
$\mathbb{R}$ when identified with vectors in $\mathbb{R}^2$.  We now show there
exist scalars $a_1,\ a_2$ such that
$$
a_1 \Re [\mu_1 S] + a_2 \Re [\mu_2 S] =S
$$
which proves the lemma.  Such scalars must solve the linear system
\begin{equation}
\begin{pmatrix} \mu_1 && \mu_2 \\
                \overline{\mu_1} && \overline{\mu_2} \end{pmatrix}
                \begin{pmatrix} a_1 \\
                                a_2 \end{pmatrix} = \begin{pmatrix} 1 \\
                                                                    0
                                                                    \end{pmatrix}
\end{equation}
Writing $\mu_1=\alpha_1 +i\beta_1$, $\mu_2=\alpha_2+i\beta_2$, a short
calculation shows that
$$
\det \begin{pmatrix} \mu_1 && \mu_2 \\
                \overline{\mu_1} && \overline{\mu_2} \end{pmatrix}
=-2i \det \begin{pmatrix} \alpha_1 && \alpha_2 \\
                          \beta_1 && \beta_2 \end{pmatrix}
$$
The latter determinant is nonzero since $\mu_1$, $\mu_2$ are linearly
independent over $\mathbb{R}$, and the system is solvable. 
$\square$

The above argument does not depend on the discreteness of $\Gamma$;
indeed it is a refinement of an argument due to J. Moorhouse
\cite{moorhouse-thesis} that the $C^*$-algebra generated by \emph{all}
M\"{o}bius transformations contains $S$.  The argument is valid for any
group which has an orbit with three accumulation points on $\mathbb{T}$;
e.g. the conclusion holds for any dense subgroup of $PSU(1,1)$.

\vskip.2in

\noindent\textbf{Proof of Theorem~\ref{T:main}}  For an automorphism $\gamma$ of
$\mathbb{D}$ set
$$
U_\gamma =(C_\gamma C_\gamma^*)^{-1/2} C_\gamma,
$$ the unitary appearing in the polar decomposition of $C_\gamma$.  By
Lemma \ref{L:factor}, $C_\gamma C_\gamma^*=T_f T_f^*$, so we may write
$$
U_\gamma = T_{|f|^{-1}} C_\gamma +K
$$
for some compact $K$.  Now by Lemma \ref{L:commute}, if $p$
is any analytic polynomial, 
\begin{align}
\nonumber U_\gamma T_p &=T_{|f|^{-1}} C_\gamma T_p +K^\prime\\
\nonumber &= T_{|f|^{-1}} T_{p\circ \gamma} C_\gamma +K^\prime\\
\nonumber &= T_{p \circ \gamma}T_{|f|^{-1}} C_\gamma +K^{\prime\prime} \\
\nonumber &= T_{p\circ \gamma}U_\gamma +K^{\prime\prime} 
\end{align}
where we have used the fact that $T_{|f|^{-1}}$
and $T_p$ commute modulo $\mathcal{K}$.  Taking adjoints and sums shows that
\begin{equation}\label{E:polycovar}
U_\gamma T_q U_\gamma^* =T_{q\circ \gamma} +K
\end{equation}
for any trigonometric polynomial $q$.  
We next show that, for M\"{o}bius transformations $\gamma, \eta$, 
$$
U_\gamma U_\eta = U_{\eta\circ\gamma} +K
$$
for some compact $K.$ To see this, write 
$$
\gamma(z)=\frac{az+b}{cz+d}, \qquad \eta(z)=\frac{ez+f}{gz+h}
$$
Then
$$
U_\gamma = T_{|cz+d|^{-1}} C_\gamma +K_1, \qquad U_\eta = T_{|gz+h|^{-1}}
C_\eta +K_2.
$$
Note that 
$$
U_{\eta\circ\gamma} = T_{|g(az+b)+h(cz+d)|^{-1}}C_{\eta\circ\gamma} +K_3.
$$
Then
\begin{align}
\nonumber U_\gamma U_\eta &= T_{|cz+d|^{-1}} C_\gamma T_{|gz+h|^{-1}} C_\eta +K\\
\nonumber     &= T_{|cz+d|^{-1}} T_{|g\gamma (z)+h|^{-1}} C_\gamma C_\eta
+K^\prime \\
\nonumber     &=T_{|g(az+b)+h(cz+d)|^{-1}}C_{\eta\circ\gamma}+K^\prime \\
\nonumber     &= U_{\eta\circ\gamma} +K^{\prime\prime}
\end{align}

Observe now that $\pg$ is equal to the C*-algebra generated by $S$ and the unitaries $U_\gamma$.  We have already shown that $\pg$ contains $S$ and each $U_\gamma$, for the converse we recall that $C_\gamma=(C_\gamma C_\gamma^*)^{1/2}U_\gamma$, and since $C_\gamma C_\gamma^*)$ lies in $C^*(S)$ by Lemma~\ref{L:factor}, we see that $C_\gamma$ lies in the C*-algebra generated by $S$ and $U_\gamma$.  Letting $\pi$ denote the quotient map $\pi:\pg\to\pg/\mathcal{K}$, it follows hat $\pg/\mathcal{K}$ is generated as a $C^*$-algebra by a copy of
$C(\mathbb{T})$ and the unitaries $\pi(U_\gamma)$, and the map
$\gamma\to\pi(U_{\gamma^{-1}})$ defines a unitary representation of
$\Gamma$.  Let $\alpha:\Gamma\to\text{Aut}(C(\mathbb{T}))$ be given by
$$
\alpha_\gamma (f)(z) = f(\gamma^{-1}(z))
$$
Then by (\ref{E:polycovar}),
\begin{equation}\label{E:cov}
\pi(U_{\gamma^{-1}})\pi(T_f)\pi(U_{\gamma^{-1}}^*)=\pi(T_{f\circ\gamma^{-1}})=\pi(T_{\alpha_\gamma(f)})
\end{equation}
for all trigonometric polynomials $f$, and hence for all $f\in
C(\mathbb{T})$ by continuity.  Thus, $\pg/\mathcal{K}$ is generated by $C(\mathbb{T})$ and
a unitary representation of $\Gamma$ satisfying the
relation~(\ref{E:cov}).  Therefore there is a surjective
$*$-homomorphism $\rho:C(\mathbb{T})\times\Gamma \to \pg/\mathcal{K}$
satisfying
$$
\rho(f)=\pi(T_f), \qquad \rho(u_\gamma)=\pi(U_{\gamma^{-1}})
$$
for all $f\in C(\mathbb{T})$ and all $\gamma\in\Gamma$.  Let
$\mathcal{J}=\text{ker}\, \rho$.  The theorem will be proved once we show
$\mathcal{J}=0$.  To do this, we use Theorems~\ref{T:as}
and~\ref{T:amenability}.  Indeed, it suffice to show that
$C(\mathbb{T})\cap\mathcal{J}=0$: then $\mathcal{J}\subset \mathcal{J}_\lambda$ by Theorem~\ref{T:as},
and $\mathcal{J}_\lambda =0$ by Theorem~\ref{T:amenability}. 

To see that $C(\mathbb{T})\cap\mathcal{J}=0$, choose $f\in C(\mathbb{T})\cap \mathcal{J}$.  Then
$\pi(T_f)=\rho(f)=0$, which means that $T_f$ is compact.  But then
$f=0$, since nonzero Toeplitz operators are non-compact.

$\square$

\subsection{Other Hilbert spaces}

In this section we prove the analogue of Theorem~\ref{T:main} for
composition operators acting on the reproducing kernel Hilbert spaces with kernel given by 
$$
k(z,w)=\frac{1}{{{(1-z\overline{w})}^n}}
$$
for integers $n\geq2$.

We let $A^2_n$ denote the Hilbert function space on $\mathbb{D}$ with kernel $k(z,w)={(1-z\overline{w})}^{-n}$.  For $n\geq2$, this space consists of those analytic functions in $\mathbb{D}$ for which
$$
\int_\mathbb{D} |f(w)|^2 (1-|w|^2)^{n-2}dA(w)
$$
is finite, and the square root of this quantity is a norm on $A^2_n$ equivalent to the norm determined by the reproducing kernel $k$ (they are equal for $n=2$).
We fix $n\geq2$ for the remainder of this section, and let $T$ denote the operator of multiplication by $z$ on $A^2_n$.  The operator  $T$ is a contractive weighted shift, and essentially normal.  Moreover, if $f$ is any bounded analytic function on $\mathbb{D}$, then multiplication by $f$ is bounded on $A^2_n$ and is denoted $M_f$.  If $\gamma$ is a M\"{o}bius transformation, then $C_\gamma$ is bounded on $A^2_n$.  We let $\pg^n$ denote the C*-algebra generated by the operators $\{C_\gamma :\gamma\in\Gamma\}$ acting on $A^2_n$.  In this section we prove:

\begin{thm}\label{T:main_gen}
The C*-algebra $\pg^n$ contains $C^*(T)$, and in particular the compact operators $\mathcal{K}(A^2_n)$.  The map defined by $\pi(C_\gamma)=|\gamma^\prime|^{n/2} u_\gamma$ extends to a $*$-homomorphism from $\pg^n$ onto $C(\partial\mathbb{D})\times\Gamma$, and there is an exact sequence
$$
0\to \mathcal{K}(A^2_n)\to \pg^n\to C(\partial\mathbb{D})\times\Gamma \to 0
$$
\end{thm}

We collect the following computations together in a single Lemma, the analogue for $A^2_n$ of the lemmas at the beginning of the previous section.
\begin{lem}\label{L:cgn_calc}For $\gamma\in\Gamma$, with $a=\gamma^{-1}(0)$ and $n$ fixed, let 
$$
f(z)=\frac{{(1-\overline{a}z)}^n}{{(1-{|a|}^2)}^{n/2}}
$$
Then $C_\gamma C_\gamma^* =M_fM_f^*$, and $M_\gamma C_\gamma = C_\gamma T$.
\end{lem}

\Prf As in the Hardy space, for the first equation it suffices to prove that
\begin{equation}\label{E:bilin_general}
\langle C_\gamma C_\gamma^* k_w, k_z\rangle = \langle T_f T_f^* k_w,
k_z\rangle
\end{equation}
for all $z,w$ in $\mathbb{D}$.  This is the same as the identity
\begin{equation}\label{E:phi_ident_general}
{\left( \frac{1}{1-\gamma(z)\overline{\gamma(w)}}\right)}^n=\frac{{(1-\overline{a}z)}^n{(1-a\overline{w})}^n}{{(1-|a|^2)}^n{(1-z\overline{w})}^n}.
\end{equation}
which is just (ref) raised to the power $n$.  As before, the second equation follows immediately from the definitions. $\square$

The proof of Theorem~\ref{T:main_gen} follows the same lines as that for the Hardy space, except that it
requires more work to show that $\pg^n$ contains $T$.  We first prove that $\pg^n$ is irreducible.
\begin{prop}\label{P:reduce}
The $C^*$-algebra $\pg^n$ is irreducible.
\end{prop}
\Prf We claim that the operator $T^n$ lies in $\pg^n$; we first show how
irreducibility follows form this and then prove the claim.  

By \cite{MR2001h:47044}, the only reducing subspaces of $T^n$ are direct sums of subspaces of the form
$$
X_k=\overline{\text{span}}\{z^{k+mn}|m=0,1,2,\dots\}
$$ for $0\leq k\leq n-1$.  Thus, since $T^n\in \pg^n$, these are the
only possible reducing subspaces for $\pg^n$.  Suppose, then, that $X$
is a nontrivial direct sum of distinct subspaces $X_k$ (i.e. $X\neq
A^2_n$).  Observe that, for each $k$, either $X$ contains $X_k$ as a summand, or $X$ is orthogonal to $X_k$.  In the latter case (which we assume holds for some $k$), the $k^{th}$ Taylor coefficient of every function in $X$ vanishes.

It is now easy to see that $X$ cannot be reducing (or even invariant) for $\pg^n$.  Indeed, if $X$ does not contain $X_0$ as a summand, then every function $f\in X$ vanishes at the origin, but if $\gamma\in \Gamma$ does not fix the origin then $C_\gamma f\notin X$.  On the other hand, if $X$ contains the scalars, then consider the operator $F(\lambda)=(1-\overline{\lambda}T)^n(1-\lambda T^*)^n$, for $\lambda\in\Lambda$.  Since $T^*$ annihilates the scalars, we have  $F(\lambda)1=(1-\overline{\lambda}z)^n=p(z)$.  But since $|\lambda|=1$, the $k^{th}$ Taylor coefficient of $p$ does not vanish, so $p\notin X$.  But since $F(\lambda)$ belongs to $\pg^n$, it follows that $X$ is not invariant for $\pg^n$ and hence $\pg^n$ is irreducible.

To prove that $T^n\in \pg^n$, we argue along the lines of
Proposition~\ref{P:shift}, though the situation is somewhat more
complicated.  Using the same reproducing kernel argument as in that Proposition, we see
that the operators
\begin{align}
F(\lambda)&=(1-\overline{\lambda}T)^n(1-\lambda T^*)^n \\
&= \sum_{j=0}^n\sum_{k=0}^n (-1)^{(j+k)} \binom{n}{j}\binom{n}{k}\overline{\lambda}^k \lambda^j T^k T^{*j} \\
&=\sum_{d=0}^n \binom{n}{d}^2 T^d T^{*d} +2\Re \sum_{m=1}^n (-1)^m\overline{\lambda}^m \sum_{m\leq k\leq n}\binom{n}{k}\binom{n}{k-m}T^k T^{*{k-m}}   \\
&= \sum_{d=0}^n \binom{n}{d}^2 T^d T^{*d} + \sum_{m=1}^n \overline{\lambda}^m E_m + \lambda^m E_m^*
\end{align}
lie in $\pg^n$ for all $\lambda$ in the limit set $\Lambda$ of $\Gamma$.  Here we have adopted the notation
$$
E_m = (-1)^m \sum_{m\leq k\leq n} \binom{n}{k}\binom{n}{k-m}T^k T^{* k-m}
$$
Note in particular that $E_n=(-1)^n T^n$. 
Forming differences as $\lambda$ ranges over $\Lambda$ shows that $\pg^n$ contains all operators of the form
$$
G(\lambda,\mu)=F(\lambda)-F(\mu) = \sum_{m=1}^n (\overline{\lambda}^m-\overline{\mu}^m)E_m +(\lambda^m-\mu^m)E_m^*
$$
for all $\lambda, \mu \in \Lambda$.  We wish to obtain $T^n$ as a
linear combination of the $G(\lambda,\mu)$; since $E_n$ is a scalar multiple of $T^n$ it suffices to  show
that there exist $2n$ pairs $(\lambda_j, \mu_j)\in \Lambda\times\Lambda$ and $2n$ scalars
$\alpha_j$ such that 
$$
E_n = \sum_{j=1}^{2n} \alpha_j G(\lambda_j, \mu_j)
$$
Let $L$ be the $2n\times 2n$ matrix whose $j^{th}$ column is
$$
\begin{pmatrix} \overline{\lambda_j}-\overline{\mu_j} \\
                \lambda_j-\mu_j \\
                \overline{\lambda_j}^2-\overline{\mu_j}^2 \\
                \lambda_j^2-\mu_j^2 \\
                \vdots \\
                \overline{\lambda_j}^n-\overline{\mu_j}^n \\
                \lambda_j^n-\mu_j^n \\
\end{pmatrix}
$$
We must therefore solve the linear system  
$$
 L \begin{pmatrix} \alpha_1 \\ 
                   \vdots \\
                   \alpha_{2n-1}\\
                   \alpha_{2n}\\
  \end{pmatrix}
 = \begin{pmatrix} 0 \\
                   \vdots \\
                   1 \\
                   0 \end{pmatrix}
$$
for which it suffices to show that the matrix $L$ is 
nonsingular for some choice of the $\lambda_j$ and $\mu_j$ in
$\Lambda$.  To do this, we fix $2n+1$ distinct points $z_0, z_1,\dots
z_{2n}$ in $\Lambda$ and set $\lambda_j=z_0$ for all $j$ and $\mu_j=z_j$
for $j=1,\dots 2n$.  The matrix $L$ then becomes the matrix whose $j^{th}$
column is 
$$
\begin{pmatrix} \overline{z_0}-\overline{z_j} \\
                z_0-z_j \\
                \overline{z_0}^2-\overline{z_j}^2 \\
                z_0^2-{z_j}^2 \\
                \vdots \\
                \overline{z_0}^n-\overline{z_j}^n \\
                z_0^n-{z_j}^n \\
\end{pmatrix}
$$
To prove that $L$ is nonsingular, we prove that its rows are linearly independent.  To this end, let $c_j$ and $d_j$, $j=1,\dots n$ be scalars such that for each $k=1, \dots 2n$,
\begin{equation}\label{E:roots}
\sum_{j=1}^n c_j (z_0^j-z_k^j)+\sum_{j=1}^n d_j (\overline{z_0}^j-\overline{z_k}^j) = 0
\end{equation}
To see that all of the $c_j$ and $d_j$ must be $0$, consider the harmonic polynomial
$$
P(z,\overline{z})=\sum_{j=1}^n c_j (z_0^j-z^j)+\sum_{j=1}^n d_j (\overline{z_0}^j-\overline{z}^j) 
$$
By (\ref{E:roots}), $P$ has $2n+1$ distinct zeroes on the unit circle, namely the points $z_0, \dots z_{2n}$.  But this means that the rational function
$$
Q(z)=\sum_{j=1}^n c_j (z_0^j-z^j)+\sum_{j=1}^n d_j (\overline{z_0}^j-\frac{1}{z^j}) 
$$
also has $2n+1$ zeroes on the circle.  But since the degree of $Q$ is at most $2n$, it follows that $Q$ must be the zero polynomial, and hence all the $c_j$ and $d_j$ are $0$.

\begin{prop}\label{P:bergshift}
$\pg^n$ contains $T$.
\end{prop}
\Prf 
We have established that for each $\lambda$ in the limit set, the operators
$$
(1-\overline{\lambda}T)^n(1-\overline{\lambda}T)^{*n}
$$
lie in $\pg^n$.  Now, since $T$ is essentially normal, the
difference
$$
(1-\overline{\lambda}T)^n(1-\overline{\lambda}T)^{*n} -
[(1-\overline{\lambda}T)(1-\overline{\lambda}T)^*]^n
$$ is compact.  Given
two positive operators which differ by a compact, their (unique)
positive $n^{th}$ roots also differ by a compact.  The positive
square root of the left-hand term above is
thus equal to 
$$
(1-\overline{\lambda}T)(1-\overline{\lambda}T)^*+K
$$ for some compact operator $K$ (depending on $\lambda$), and lies in $\pg^n$.  Forming
linear combinations as in Proposition~\ref{P:shift}, we conclude that
$T+K$ lies in $\pg^n$ for some compact $K$.  Now, as $\pg^n$ is
irreducible and contains a nonzero compact operator (namely, the
self-commutator of $T+K$), it contains all the compacts, and hence
$T$. $\square$

We now prove the analogue of Theorem~\ref{T:main} for $\pg^n$:

\begin{thm}\label{T:berg_main}  
  Let $\Gamma$ be a non-elementary Fuchsian group, and let $\pg^n$
  denote the $C^*$-algebra generated by the set of composition
  operators on $A^2$ with symbols in $\Gamma$.  Then there is an exact
  sequence
\begin{equation}\begin{CD}
0 @>>> \mathcal{K} @>\iota>> \pg^n @>\pi >> C(\mathbb{T})\times \Gamma @>>> 0
\end{CD}
\end{equation}

\end{thm}
\Prf The proof is similar to that of Theorem~\ref{T:main}; we define
the unitary operators $U_\gamma$ in the same way (using the polar
decomposition of $C_\gamma$), and check that he map that sends $\gamma$ to
$U_{\gamma^{-1}}$ is a unitary representation of $\Gamma$ on $A^2_n$,
modulo $\mathcal{K}$.  The computation to check the covariance
condition is essentially the same.  As $\pg^n$ is generated by the
$T$ and the unitaries $U_\gamma$, the quotient $\pg^n/\mathcal{K}$ is
generated by a copy of $C(\partial\mathbb{D})$ (since
$\partial\mathbb{D}$ is the essential spectrum of $T$)
and a representation of $\Gamma$ which satisfy the covariance
condition; the rest of the proof proceeds exactly as for $H^2$, up until the final step.  At the final step in that proof, we have a Toeplitz operator $T_f$ which is compact; this implies that $f=0$ in the $H^2$ case but not for the spaces $A^2_n$.  However the symbol of a compact Toeplitz operator on $A^2_n$ must vanish at the boundary of $\mathbb{D}$, so we may still conclude that $f=0$ on $\partial\mathbb{D}$ and the proof is complete.$\square$


\section{$K$-Homology of $C(\partial\mathbb{D})\times\Gamma$}

\subsection{$Ext$ classes of $\pg$ and $\pg^n$}
In this section we prove that the extensions of $C(\partial\mathbb{D})\times
\Gamma$ determined by $\pg$ and $\pg^n$ determine the same element of
the group $\text{Ext}(C(\partial\mathbb{D})\times\Gamma)$, and we prove that
$\pg$ and $\pg^n$ are isomorphic as $C^*$-algebras.  (In general neither
of these statements implies the other.)  We borrow several ideas from
[R\o rdam] concerning the classification of extensions by the
associated six-term exact sequence in $K$-theory.  The reader is
referred to \cite{MR2001g:46001} for an introduction to $K$-theory and
\cite[Chapter 15]{MR88g:46082} for the basic definitions
and theorems concerning extensions of $C^*$-algebras and the $Ext$ group.

\begin{thm}\label{T:extclasses}
The extensions of $C(\partial\mathbb{D})\times\Gamma$ determined by $\pg$ and
$\pg^n$ define the same element of the group
$\text{Ext}(C(\partial\mathbb{D})\times\Gamma)$.
\end{thm}

We will show that the Busby maps of the extensions determined by $\pg$
and $\pg^n$ are strongly unitarily equivalent.  To do this we introduce an
integral operator $V\in\mathcal{B}(A^2_n,A^2_{n-1})$
and prove the following lemma, which describes the properties of $V$
needed in the proof of the theorem.
\begin{lem}\label{L:v_operator}Consider the integral operator defined on analytic polynomials $f(w)$ by
$$
(Vf)(z)=\int_{\mathbb{D}}\frac{f(w)}{(1-\overline{w}z)}(1-|w|^2)^{-1/2}\, dA(w)
$$
where $dA$ is normalized Lebesgue area measure on $\mathbb{D}$.  
Then:
\begin{enumerate}
\item For each $n\geq0$, $V$ extends to a bounded diagonal operator from $A^2_{n+1}$ to $A^2_n$, and there exists a compact operator $K_n$ such that $V+K_n$ is unitary.
\item If $T_{n+1}$ and $T_n$ denote multiplication by $z$ on $A^2_{n+1}$ and $A^2_n$ respectively, then $VT_{n+1}-T_nV$ is compact. 
\item If $g$ is continuous on $\overline{\mathbb{D}}$ and analytic in $\mathbb{D}$, then the integral operator 
$$
(V_{\overline{g}}f)(z)=\int_\mathbb{D}\frac{f(w)\overline{g(w)}}{1-\overline{w}z}(1-|w|^2)^{-1/2}dA(w)
$$
is bounded from $A^2_{n+1}$ to $A^2_n$, and is a compact perturbation of $VM_g^*$.
\end{enumerate}
\end{lem}

\Prf If $f(w)=w^k$, then the integrand is absolutely convergent for each $z\in\mathbb{D}$, and we calculate
\begin{align}
(Vf)(z) &=\int_\mathbb{D}\frac{w^k}{1-\overline{w}z}(1-|w|^2)^{-1/2}dA(w) \\
&= \sum_{j=0}^\infty\int_\mathbb{D}w^k\overline{w}^j z^j (1-|w|^2)^{-1/2}dA(w)\\
&= \left(\int_\mathbb{D}|w|^{2k} (1-|w|^2)^{-1/2}dA(w)\right) z^k \\
&= \alpha_k z^k 
\end{align}
It is known that the sequence $(\alpha_k)$ defined by the above integral is asymptotically $(k+1)^{-1/2}$.  Each space $A^2_n$ has an orthonormal basis of the form $\beta_k z^k$, where the sequence $\beta_k$ is asymptotically $(k+1)^{(n-1)/2}$.  It follows that $V$ intertwines orthonormal bases for $A^2_{n+1}$ and $A^2_n$ modulo a compact diagonal operator, which establishes the first statement.  Moreover, the second statement also follows, by observing that the operators $T_{n+1}$ and $T_n$ are weighted shifts with weight sequences asymptotic to $1$.  

To prove the last statement, first suppose $g(w)=w$, and let $\beta_k z^k$ be an orthonormal basis for $A^2_{n+1}$.  Then a direct computation shows that $V_{\overline{g}}$ is a weighted backward shift with weight sequence $\alpha_{j+1} \beta_j /\beta_{j+1}$, $j=0, 1, \dots$.  Thus with respect to this basis, $V^*V_{\overline{z}}$ is a weighted shift with weight sequence $\alpha_{j+1}\beta_j /\alpha_j\beta_{j+1}$.  Since $\lim_{j \to\infty}\alpha_{j+1}/\alpha_j=1$, it follows that $V^*V_{\overline{z}}$ is a compact perturbation of $M_z^*$.  By linearity, the lemma holds for polynomial $g$.  Finally, if $p_n$ is a sequence of polynomials such that $M_{p_n}$ converges to $M_g$ in the operator norm (which is equal to the supremum norm of the symbol), then $M_{p_n}$ converges to $M_g$ in the essential norm, and $p_n$ converges to $g$ uniformly, so that $V_{\overline{p_n}}\to V_{\overline{g}}$ in norm and the result follows.

\noindent\textbf{Proof of Theorem~\ref{T:extclasses}}

We will prove that the Busby maps of the extensions determined by
$\pg$ and $\pg^n$ are strongly unitarily equivalent.  It suffices to
prove, for each fixed $n$, the equivalence between $\pg^n$ and
$\pg^{n+1}$.  We first establish some notation: for a function $g$ in
the disk algebra, we let $M_g$ denote the multiplication operator with
symbol $g$ acting on $A^2_{n+1}$, and $\widetilde{M}_g$ the
corresponding operator acting on $A^2_n$.  Similarly, for the
composition operators $C_\gamma$ and the unitaries $U_\gamma$ on
$A^2_{n+1}$, a tilde denotes the corresponding operator on $A^2_n$.
Now, since each of the algebras $\pg^n, \pg^{n+1}$ is generated by the
operators $M_z$ (resp. $\widetilde{M}_z$) and the unitaries $U_\gamma$
(resp. $\widetilde{U}_\gamma$) on the respective Hilbert spaces, it
suffices to produce a unitary $U$ (or indeed a compact perturbation of
a unitary) from $A^2_{n+1}$ to $A^2_n$ such that the operators
$UM_z-\widetilde{M}_zU$ and $UU_\gamma-\widetilde{U}_\gamma U$ are
compact for all $\gamma\in\Gamma$.  In fact we will prove that the
operator $V$ of the previous lemma does the job.

To prove this, we will first calculate the operator $\widetilde{C}_\gamma VC_{\gamma^{-1}}$.
Written as an integral operator, 
$$
(\widetilde{C}_\gamma VC_{\gamma^{-1}}
f)(z)=\int_\mathbb{D}\frac{f(\gamma^{-1}(w))}{1-\overline{w}\gamma(z)}(1-|w|^2)^{-1/2}\,dA(w)
$$
Applying the change of variables $w\to\gamma(w)$ to this integral, we
obtain
$$
(\widetilde{C}_\gamma VC_{\gamma^{-1}}
f)(z)=\int_\mathbb{D}\frac{f(w)}{1-\overline{\gamma(w)}\gamma(z)}(1-|\gamma(w)|^2)^{-1/2}|\gamma^\prime(w)|^2\,dA(w)
$$
A little algebra shows that 
$$
1-|\gamma(w)|^2 = (1-|w|^2)|\gamma^\prime(w)|
$$
Furthermore, multiplying and dividing the identity (\ref{E:phi_ident})
by $1-\overline{a}w$, we obtain
$$
\frac{1}{1-\overline{\gamma(w)}\gamma(z)}=\frac{1-\overline{a}z}{(1-\overline{a}w)(1-\overline{w}z)}|\gamma^\prime(w)|^{-1}
$$
Thus the integral may be transformed into 
$$
(\widetilde{C}_\gamma VC_{\gamma^{-1}} f)(z)=\int_\mathbb{D}\frac{f(w)}{1-\overline{w}z}\frac{1-\overline{a}z}{1-\overline{a}w}|\gamma^\prime(w)|^{1/2}(1-|w|^2)^{-1/2}\,dA(w)
$$ 
Since $\gamma^\prime$ is analytic and non-vanishing in a neighborhood of the closed disk, there exists a function $\psi$ in the disk algebra such that $|\psi|^2=|\gamma^\prime|^{1/2}$ (choose $\psi$ to be a branch of $(\gamma^\prime)^{1/4}$).  Thus we may rewrite this integral as
\begin{align}
(\widetilde{C}_\gamma VC_{\gamma^{-1}} f)(z) &=\int_\mathbb{D}\frac{f(w)}{1-\overline{w}z}\frac{1-\overline{a}z}{1-\overline{a}w}\psi(w)\overline{\psi(w)}(1-|w|^2)^{-1/2}\,dA(w)\\
&= M_{1-\overline{a}z}V_{\overline{\psi}}M_{(1-\overline{a}z)^{-1}\psi(z)}
\end{align} 
Now, applying the lemma to $V_{\overline{\psi}}$ and using the fact that $V$ intertwines the multiplier algebras of $A^2_{n+1}$ and $A^2_n$ modulo the compacts, we have
\begin{align}
\widetilde{C}_\gamma VC_{\gamma^{-1}} &= M_{1-\overline{a}z}V_{\overline{\psi}}M_{(1-\overline{a}z)^{-1}\psi(z)}\\
&= M_{1-\overline{a}z} VM_\psi^*M_\psi M_{(1-\overline{a}z)^{-1}} +K\\
&= VM_\psi^*M_\psi +K
\end{align}
modulo the compacts.  Multiplying on the right by $C_{\gamma^{-1}}$ and on the left by $V^*$ we get
\begin{equation}\label{E:VCV}
V^*\widetilde{C}_\gamma V =M_\psi^* M_\psi C_\gamma +K
\end{equation}
Since $\psi^2(z)=(1-|a|^2)^{1/2}(1-\overline{a}z)^{-1}$, the calculations in the proof of~\ref{T:main} show  that
$$
(C_\gamma C_\gamma^*)^{-1/2}=M_{\psi^{n+1}}M_{\psi^{n+1}}^* \quad\text{and}\quad (\widetilde{C}_\gamma \widetilde{C}_\gamma^*)^{-1/2}=\widetilde{M}_{\psi^{n}}\widetilde{M}_{\psi^{n}}^*
$$
Using the fact that $V$ intertwines the C*-algebras generated by $M_z$ and $\widetilde{M}_z$ modulo $\mathcal{K}$, and that these algebras are commutative modulo $\mathcal{K}$, we conclude after multiplying both sides of~(\ref{E:VCV}) on the left by $M_{\psi^n}M_{\psi^n}^*$ that
$$
V^*\widetilde{U}_\gamma V = V^*\widetilde{M}_{\psi^{n}}\widetilde{M}_{\psi^{n}}^* \widetilde{C}_\gamma V+K=M_{\psi^{n+1}}M_{\psi^{n+1}}^* C_\gamma +K =U_\gamma +K
$$
which completes the proof.

\subsection{$C^*$-algebras isomorphic to $\pg$}

In this section we prove that all unital extensions of
$C(\partial\mathbb{D})\times\Gamma$ which determine the same
\emph{Ext} class as $\mathcal{C}_\Gamma$ have isomorphic pull-backs,
i.e. the middle term in the exact sequence is isomorphic to
$\mathcal{C}_\Gamma$.  TO prove this, we apply Voiculescu's theorem to
obtain a stable isomorphism and use an analysis of the $K_0$-classes
of finite projections in $\pg$ to obtain a unital isomorphism.

We recall that if $A$ is a unital, nuclear C*-algebra, an extension of $A$ by $\mathcal{K}$ is called \emph{unital} if the Busby map $\tau:A\to \mathcal{Q}(\mathcal{H})$ is unital.   A trivial extension is called \emph{strongly unital} $\tau$ lifts to a unital homomorphism $\rho:A\to \mathcal{B}(\mathcal{H})$.  (Not all unital trivial extensions are strongly unital.)  We begin with the following lemma, which is a special case of \cite[Proposition
5.1]{MR99b:46108}.
\begin{lem}
  Let $A$ be a unital $C^*$-algebra, and let $\tau_0=\pi\circ \alpha$,
  $\tau_1:A\to \mathcal{Q}(\mathcal{H})$ be unital extensions, with $\tau_0$
  trivial.  Then there is an isometry $v\in \mathcal{B}(\mathcal{H})$ such that
  $\alpha(1_A)=vv^*$.  Set $\alpha^\prime(a)=v^*\alpha(a)v$ and set
  $\tau^\prime_0 = \pi\circ \alpha^\prime$, so that $\tau_0^\prime$ is
  a strongly unital trivial extension.  Let
\begin{equation}\nonumber\begin{CD}
e:\quad  0 @>>> \mathcal{K} @>>> E @>\psi >> A @>>> 0
\end{CD}
\end{equation}
\begin{equation}\nonumber\begin{CD}
e^\prime:\quad  0 @>>> \mathcal{K} @>>> E^\prime @>\psi^\prime >> A @>>> 0
\end{CD}
\end{equation}
be the extensions with Busby maps $\tau_1 \oplus \tau_0$,
respectively, $\tau_1 \oplus \tau_0^\prime$.  Let $s_1, s_2$ be
isometries in $\mathcal{B}(\mathcal{H})$ with $s_1s_1^* + s_2s_2^* =1$ and set
$w=s_1s_1^* +s_2 vs_2^*$, $p=ww^*$.  Then $\beta(b)=w^*bw$ and
$\eta(x)=w^*xw$ define an isomorphism
\begin{equation}\nonumber\begin{CD}
0 @>>> p\mathcal{K} p @>>> pEp @>\psi >> A @>>> 0 \\
& & @V\beta VV @VV\eta V @|   \\
0 @>>> \mathcal{K} @>>> E^\prime @>>\psi^\prime > A @>>> 0 
\end{CD}
\end{equation}
\end{lem}

\Prf \cite{MR99b:46108}
\begin{thm}\label{T:classify}
Let $x\in \rm{Ext}(C(\partial\mathbb{D})\times\Gamma)$ denote the class of the extension
\begin{equation}\nonumber\begin{CD}
0 @>>> \mathcal{K} @>>> \pg @>>> C(\partial\mathbb{D})\times\Gamma @>>> 0.
\end{CD}
\end{equation}
If $e\in \rm{Ext}(C(\partial\mathbb{D})\times\Gamma)$ is a unital extension represented by 
\begin{equation}\nonumber\begin{CD}
0 @>>> \mathcal{K} @>>> E  @>>> C(\partial\mathbb{D})\times\Gamma @>>> 0
\end{CD}
\end{equation}
and $e=x$, then $E\cong \pg$ as $C^*$-algebras.
\end{thm}
\Prf We first prove that the index homomorphism $\partial$ from
$K_1(C(\partial\mathbb{D})\times\Gamma)$ to $\mathbb{Z}$ determined by the extension $x$ is
surjective (and hence so is the homomorphism determined by $e$).
First, consider the function $f(z)=z$; this function is a unitary in
$C(\partial\mathbb{D})\times\Gamma$ lying in the canonical subalgebra
isomorphic to $C(\partial\mathbb{D})$.  By construction, this
unitary lifts to the unilateral shift $S$ on $H^2$, which is a
Fredholm isometry of index $-1$.  Thus $\partial[z]_1 =-1$, and since
$-1$ generates $\mathbb{Z}$, $\partial$ is surjective.

Since $\partial$ is surjective, the six-term exact sequence in $K$-theory
associated to $E$
\begin{equation}\begin{CD}
&0 @>>> &K_1(E) @>>> &K_1(C(\partial\mathbb{D})\times\Gamma) \\
&@AAA & & & & @VV\partial V \\
&K_0(C(\partial\mathbb{D})\times\Gamma) @<<< &K_0(E) @<<< &\mathbb{Z}  
\end{CD}
\end{equation}
splits into the two sequences
\begin{equation}\begin{CD}
0 @>>> K_1(E) @>>> K_1(C(\partial\mathbb{D})\times\Gamma) @>\partial >> \mathbb{Z} @>>> 0
\end{CD}
\end{equation}
\begin{equation}\begin{CD}
0 @>>> K_0(E) @>>> K_0(C(\partial\mathbb{D})\times\Gamma) @>>>  0
\end{CD}
\end{equation}
Thus $K_0(E)\cong K_0(C(\partial\mathbb{D})\times\Gamma)$, and as the isomorphism is induced
by the quotient map which annihilates the compacts, it follows that
$[p]_0=0$ for any finite projection $p\in E$.

Now let $\tau_1$ and $\tau_2$ denote the Busby maps associated to
$\pg$ and $E$ respectively.  Since $\pg$ and $E$ determine the same
element in $\text{Ext}(C(\partial\mathbb{D})\times\Gamma)$, there exist unital trivial
extensions $\sigma_1$ and $\sigma_2$ such that $\tau_1 \oplus \sigma_1
\sim_u \tau_2 \oplus \sigma_2$. Let $E_j$ denote the extension with
Busby map $\tau_j \oplus \sigma_j$; we have $E_1\cong E_2$.  Using
$\sigma_1$ and $\sigma_2$, we may choose isometries $v_1$ and $v_2$ in
$\mathcal{B}(\mathcal{H})$ as in the lemma to obtain strongly unital extensions
$\sigma_1^\prime$ and $\sigma_2^\prime$.  Let $E_j^\prime$ be the
extension of $C(\partial\mathbb{D})\times\Gamma$ with Busby map $\tau_j \oplus
\sigma_j^\prime$, $j=1,2$.  By Voiculescu's theorem, $\tau_j \sim_u
\tau_j \oplus \sigma_j^\prime$, and it follows (since unitarily
equivalent Busby maps determine isomorphic extensions) that
$E_1^\prime \cong \pg$ and $E_2^\prime \cong E$.  Applying the lemma
to $E_1$ and $E_2$, we obtain projections $p_j\in E_j$ such that $p_j
E_j p_j \cong E_j^\prime$.  We claim that $[p_j]_0 = [1_{E_j}]_0$;
from this the theorem follows, since we then have $p_j E_j p_j \cong
E_j$.

To prove the claim, observe that for the isometries $v_j$, the
projections $1-v_jv_j^*$ are finite, and hence so are the projections
$1-p_j$:
\begin{align}\nonumber
1-p = 1-ww^* &= 1-[s_1s_1^* + s_2vv^*s_2^*] \\
\nonumber             &= s_2s_2^* - s_2vv^*s_2^* \\
\nonumber             &= s_2[1-vv^*]s_2^* 
\end{align}
which is finite.  Thus $1-p_j$ lies in $E_j$ (since $E_j$ contains the
compacts) and $p_j\in Ej$ since $E_j$ is unital. Finally, $[1-p_j]_0
=0$, and $[p_j]_0 =[1_{E_j}]_0$.  $\square$

\subsection{Emerson's construction}
When $\Gamma$ is cocompact, there is another extension of
$C(\partial\mathbb{D})\times\Gamma$ which was constructed by Emerson
\cite{emerson-thesis}, motivated by work of Kaminker and Putnam
\cite{MR98f:46056}; we will show that the $Ext$-class of $\pg$ differs
from the $Ext$-class of this extension.

The extension is constructed as follows: since $\Gamma$ is cocompact,
we may identify $\mathbb{T}$ with the Gromov boundary $\partial\Gamma$.  Let
$f$ be a continuous function on $\partial\Gamma$, extend it
arbitrarily to $\Gamma$ by the Tietze extension theorem, denote this
extended function by $\tilde{f}$.  Let $e_x, x\in\Gamma$ be the
standard orthonormal basis for $\ell^2(\Gamma)$, and let $u_\gamma$,
$\gamma\in\Gamma$ denote the unitary operator of left translation on
$\Gamma$.  Define a map $\tau:C(\partial\mathbb{D})\times\Gamma\to \mathcal{Q}(\ell^2(\Gamma))$ by
$$
\tau(f)e_x = \tilde{f}(x)e_x
$$
and
$$
\tau(\gamma)e_x =u_\gamma e_x = e_{\gamma x},
$$
it can be shown that these expressions are well defined modulo the
compact operators and determine a $*$-homomorphism from $C(\partial\mathbb{D})\times\Gamma$ to
the Calkin algebra of $\ell^2(\Gamma)$, i.e. an extension of
$C(\partial\mathbb{D})\times\Gamma$ by the compacts.  Let $\pi$ denote the quotient map
$\pi:\mathcal{B}(\ell^2(\Gamma))\to \mathcal{Q}(\ell^2(\Gamma))$, and consider the
pull-back $C^*$-algebra $E$:
\begin{equation}\begin{CD}
E @>\tilde{\tau}>> C(\partial\mathbb{D})\times\Gamma \\
@VVV   @VV\tau V \\
\mathcal{B}(\ell^2(\Gamma)) @>>\pi> \mathcal{Q}(\ell^2(\Gamma))
\end{CD}
\end{equation}

We will show this extension is distinct from $\pg$ by showing that it
induces a different homomorphism from $K_1(C(\partial\mathbb{D})\times\Gamma)$ to $\mathbb{Z}$.
Indeed, consider the class $[z]_1$ of the unitary $f(z)=z$ in
$K_1(C(\partial\mathbb{D})\times\Gamma)$; the extension belonging to
$\pg$ sends this class to $-1$.  On the other hand, for the extension
$\tau$ we claim the function $z$ lifts to a unitary in $E$; it follows
that $\partial ([z]_1) =0$ in this case.  To verify the claim, note
that $z$ lifts to a diagonal operator on $\ell^2(\Gamma)$, and if
$(x_n)$ is a subsequence in $\Gamma$ tending to the boundary point
$\lambda\in\mathbb{T}$, then $\tilde{f}(x_n)\to f(\lambda)=\lambda$.  We may
thus choose all the $\tilde{f}(x)$ nonzero, and by dividing
$\tilde{f}(x)$ by its modulus, we may assume they are all unimodular.
Thus $z$ lifts to the unitary $diag(\tilde{f}(x))\oplus z$ in $E$.

\subsection{Lott's construction}

In this section we relate the extension $\tau$ given by
$\mathcal{C}_\Gamma$ to an extension recently constructed by J. Lott \cite{lott-preprint}.  In particular it will follow that (up to tensoring with
$\mathbb{Q}$) $\tau$ lies in the range of the Baum-Douglas-Taylor
boundary map
$$
\partial: K^0 (C(\overline{\mathbb{D}})\times\Gamma, C(\partial\mathbb{D})\times\Gamma )\to K^1(C(\partial\mathbb{D})\times\Gamma)
$$
We first describe a construction of the extension $\sigma_+$ of
\cite{lott-preprint}.  Let $\mathcal{D}$ denote the Hilbert space of
analytic functions on $\mathbb{D}$ with finite Dirichlet integral
$$
D(f)=\int_\mathbb{D}|f^\prime(z)|^2dA(z)
$$
equipped with the norm
$$
\|f\|^2=|f(0)|^2 + \int_\mathbb{D}|f^\prime(z)|^2dA(z)=|f(0)|^2+D(f)
$$
If $f$ is represented in $\mathbb{D}$ by the Taylor series $\sum_{n=0}^\infty a_n z^n$, then this norm is given by
$$
\|f\|^2 = |a_0|^2 + \sum_{n=1}^\infty n |a_n|^2
$$
The operator of multiplication by $z$ on $\mathcal{D}$, denoted $M_z$,  is a weighted shift with weight sequence asymptotic to $1$, and hence is unitarily equivalent to a  compact perturbation of the unilateral shift on $H^2$.  It follows that there is a *-homomorphism $\rho:C(\partial\mathbb{D})\to \mathcal{Q}(\mathcal{D})$ with $\rho(z)=\pi(M_z)$.  Now, by changing variables one checks that if $\gamma$ is a M\"{o}bius transformation, $D(f\circ\gamma)=D(f)$.  Let $\mathcal{D}_0$ denote the subspace of $\mathcal{D}$ consisting of those functions which vanish at the origin.  It then follows from the definition of the norm in $\mathcal{D}$ that the operators 
$$
u_\gamma(f)(z)=f(\gamma(z))-f(\gamma(0))
$$ 
are unitary on $\mathcal{D}_0$, and form a unitary representation of $\Gamma$.  We extend this representation to all of $\mathcal{D}$ by letting $\Gamma$ act trivially on the scalars.  Moreover, it is simple to verify (by noting that $u_\gamma$ is a compact perturbation of the composition operator $C_\gamma$) that for all $\gamma\in\Gamma$, 
$$
u_\gamma M_z u_\gamma ^* = M_{\gamma(z)}
$$
modulo compact operators.  Arguing as in the proof of Theorem~\ref{T:main}, we conclude that the pair $(\rho(f), \pi(u_\gamma))$ determines an injective *-homomorphism from $C(\partial\mathbb{D})\times\Gamma$ to $\mathcal{Q}(\mathcal{D})$, which is unitarily equivalent to the Busby map $\sigma_+$ of \cite{lott-preprint}.  

We may now state the main theorem of this subsection:

\begin{thm}\label{T:lott}
The Busby maps $\tau$ and $\sigma_+$ are unitarily equivalent.
\end{thm}

\Prf  We first show that the Busby map $\tau: \mathcal{C}_\Gamma\to \mathcal{Q}(H^2)$ lifts to a completely positive map $\eta:C(\partial\mathbb{D})\times\Gamma \to\mathcal{B}(H^2)$.  Define a unitary representation of $\Gamma$ on $L^2(\partial\mathbb{D})$ by
$$
U(\gamma^{-1})=M_{|\gamma^\prime|^{1/2}}C_\gamma
$$ Together with the usual representation of $C(\partial\mathbb{D})$
as multiplication operators on $L^2$, we obtain a covariant
representation of $(\Gamma, \partial\mathbb{D})$ which in turn
determines a representation $\rho:C(\partial\mathbb{D})\times\Gamma
\to \mathcal{B}(L^2)$.  Letting $P$ denote the Riesz projection $P:L^2\to H^2$, we next claim that the commutator
$[\rho(a),P]$ is compact for all $a\in
C(\partial\mathbb{D})\times\Gamma$, and so the pair $(\rho, P)$ is an
abstract Toeplitz extension of $C(\partial\mathbb{D})\times\Gamma$ by
$\mathcal{K}$.  To see this, it suffices to prove the compactness of the commutators
$$
[\pi(f), P] \quad \text{and} \quad [U(\gamma),P].
$$ 
Now, it is well known that $[M_f, P]$ is compact, and as $U(\gamma)$
has the form $M_g C_{\gamma^{-1}}$ it suffices to check that $[C_\gamma,
P]$ is compact.  It is easily checked that this latter commutator is
rank one.  Indeed, the range of $P$ is invariant for $C_\gamma$, so
$[C_\gamma, P]= PC_\gamma P^\bot$.  If we expand $h\in
L^2(\partial\mathbb{D})$ in a Fourier series
$$
h \sim \sum_{n\in\mathbb{Z}} \hat{h}(n)e^{in\theta}
$$ 
then a short calculation shows that 
$$
PC_\gamma P^\bot h \sim (\sum_{n<0} \hat{h}(n)\overline{\gamma(0)}^{|n|})
\cdot 1 
$$ 
so $PC_\gamma P^\bot$ is rank one.

Thus, we have established that the completely positive map $\eta:C(\partial\mathbb{D})\times\Gamma \to\mathcal{B}(H^2)$ given by
$$
\eta(a)=P\rho(a)P
$$
is a homomorphism modulo compacts, and the calculations in the proof of Theorem~\ref{T:main} show that the map
$$
\tau(a)=\pi(P\rho(a)P)
$$
coincides with the Busby map associated to $\pg$.  Thus, $\eta$ is a completely positive lifting of $\tau$ as claimed.

With the lifting $\eta$ in hand, to prove the unitary equivalence of
$\tau$ and $\sigma_+$ it will therefore suffice to exhibit an operator
$V:H^2\to \mathcal{D}$ such that $V$ is a compact perturbation of a
unitary, and such that
$$
\pi(V\eta(a)V^*)=\sigma_+(a)
$$ for all $a\in C(\partial\mathbb{D})\times\Gamma$.  Since the
crossed product is generated by the function $f(z)=z$ and the formal
symbols $[\gamma]$, it suffices to establish the above equality on
these generators.  In fact, we may use the operator $V$ of
Lemma~\ref{L:v_operator}.  Since the map $z^n \to n^{-1/2}z^n$ is
essentially unitary from $H^2$ to $\mathcal{D}$, the proofs of
statements (1) and (2) of Lemma~\ref{L:v_operator} are still valid.
The conclusion of statement (3) holds provided the hypothesis on $g$
is strengthened, by requiring that $g$ by analytic in a neighborhood
of $\overline{\mathbb{D}}$.  Moreover, the arguments of
Theorem~\ref{T:extclasses} still apply, since the proof applies
statement (3) of Lemma~\ref{L:v_operator} only to the function $\psi$,
which is indeed analytic across the boundary of $\mathbb{D}$.  Thus,
the arguments in the proof of Theorem~\ref{T:extclasses} prove that
$V$ intertwines (modulo compacts) multiplication by $z$ on $H^2$ and
$\mathcal{H}$, and also intertwines (modulo compacts) $C_\gamma$
acting on $\mathcal{H}$ with $U_\gamma$ on $H^2$.  Since $u_\gamma$ is
a compact perturbation of $C_\gamma$ on $\mathcal{H}$, it follows that
$V$ intertwines $U_\gamma$ and $u_\gamma$.

Now, the Busby map $\sigma_+$ takes the function $z$ to the image of
$M_z$ in the Calkin algebra $\mathcal{Q}(\mathcal{H})$.  Since
$\eta(z)=M_z\in\mathcal{B}(H^2)$, the intertwining property of $V$
(modulo compacts) may be written as
$$
\pi(V\eta(z)V^*)=\sigma_+(z)
$$ Similarly, since $\eta$ applied to the formal symbol $[\gamma]$
(viewed as a generator of $C(\partial\mathbb{D})\times\Gamma$) is
$U_\gamma$, the intertwining property for $V$ with respect to
$U_\gamma$ and $u_\gamma$ reads
$$
\pi(V\eta([\gamma])V^*)=\sigma_+([\gamma])
$$ 
Thus the equivalence of $\tau$ and $\sigma_+$ on generators is
established, which proves the theorem. $\square$

We conclude by observing that the covariant representation on $L^2$
described in the proof of the previous theorem gives rise to an
equivariant $KK_1$-cycle for $C(\partial\mathbb{D})$.  Indeed, such a
cycle consists of a triple $(U, \pi, F)$ where $(U, \pi)$ is
a covariant representation on a Hilbert space $\mathcal{H}$ and
$F$ is a bounded operator on $\mathcal{H}$ such that the
operators 
$$
F^2-I,\  F-F^*,\  [U(\gamma),F], \text{ and } [\pi(f),F]
$$
are compact for all $\gamma\in\Gamma$ and
$f\in C(\partial\mathbb{D})$.  The computations in the previous proof
show that the triple $(U, \pi, 2P-I)$ satisfies all of these
conditions, and essentially the same unitary equivalence argument
shows that this cycle represents (up to a scalar multiple) the class
of \cite[Section 9.1]{lott-preprint}.

\bibliographystyle{plain} 
\bibliography{compext} 

\end{document}